\documentclass[12pt]{article}
\usepackage{amsmath}
\usepackage{latexsym}
\usepackage{amssymb}
\usepackage{a4wide}
\newtheorem{Defn}{Definition}
\newtheorem{Remark}{Remark}
\newtheorem{Note}{Note}
\newtheorem{prop}{Proposition}

\newtheorem{cor}{Corollary}
\newtheorem{Example}{Example}
\newtheorem{Examples}{Examples}
\newtheorem{Problems}{Problems}

\newtheorem{Problem}{Problem}
\newtheorem{Number}{\!\!}

\newenvironment{proof}{{\noindent\bf Proof.}}%
                  {\nopagebreak\hspace*{\fill}$\Box$\medskip\medskip\par}   
\newcommand{\Punkt}{\nopagebreak\hspace*{\fill}$\Box$}

\newcommand{\tensor}{\otimes}

\newcommand{\n}{\rm}

\newcommand{\mto}{\mapsto}

\newcommand{\N}{{\mathbb N}}

\newcommand{\R}{{\mathbb R}}

\newcommand{\Q}{{\mathbb Q}}

\newcommand{\C}{{\mathbb C}}

\newcommand{\K}{{\mathbb K}}

\newcommand{\one}{\mbox{\rm \bf 1}}

\newcommand{\sub}{\subseteq}

\newcommand{\Repart}{\mbox{\n Re}}
\newcommand{\Impart}{\mbox{\n Im}}

\newcommand{\repart}{\mbox{\n\footnotesize Re}}

\begin{document}
%
%
\begin{center}
{\Large \bf Examples of differentiable mappings\vspace{2mm}
into\\
non-locally convex spaces}\vspace{4.1 mm}\\
{\bf Helge Gl\"{o}ckner}\vspace{.7mm}
\end{center}
\noindent{\bf Abstract.\/}
Examples of differentiable mappings into real or
complex topological vector spaces
with specific
properties are given, which illustrate
the differences between differential calculus
in the locally
convex and the non-locally convex
case.\footnote{Classification: 58C20 (main); 26E20, 46A16, 46G20.}\\[3mm]
%
%
%
%
{\bf Introduction.$\;$}
Beyond the familiar theories
of differentiation in real and complex\linebreak
locally convex spaces
(\cite{Kel}, \cite{KaM}),
a comprehensive theory of $C^r_\K$-maps
between open subsets of topological vector spaces
over arbitrary non-discrete topological fields~$\K$
has recently been developed~\cite{Ber}.
%
Surprisingly large parts of classical differential
calculus remain intact for these maps.
For example,
being $C^r_\K$ is a local property,
the Chain Rule holds,
and $C^r_\K$-maps admit
finite order Taylor expansions~\cite{Ber}.
Furthermore, when $\K$ is a complete valued field,
implicit function theorems for $C^r_\K$-maps from topological
$\K$-vector spaces to Banach spaces are available~\cite{IMP}.
All basic constructions of infinite-dimensional Lie theory
(linear Lie groups, mapping groups,
diffeomorphism groups) work just as well
over general topological fields, valued fields,
or at least local fields~\cite{ARB}.\\[3mm]
In the real locally convex case, the $C^r_\R$-maps are precisely
the $C^r$-maps in the sense of Michal and Bastiani
(also known as Keller's $C^r_c$-maps~\cite{Kel}).
The Fundamental Theorem of Calculus
holds for such maps; in particular,
mappings whose differentials vanish at each point
have to be locally constant.
A mapping into a complex locally convex
space is of class $C^\infty_\C$ if and only if it
is complex analytic in the usual sense
(as in~\cite{BaS}). Thus, every $C^\infty_\C$-map
into a locally convex space is given locally
by its Taylor series, and the Identity Theorem
holds for such maps. Furthermore,
it is known that every $C^1_\C$-map
into a complete complex locally
convex space is automatically of class $C^\infty_\C$
(see \cite{Ber} for all of this).\\[3mm]
The purpose of this note is to describe examples
showing that these facts
become false
for mappings with non-locally convex ranges.
Thus, for suitable non-locally convex topological vector spaces
$E$, we encounter a smooth injection
$\R\to E$ whose derivative vanishes identically;
we present a $C^\infty_\C$-map $\C\to E$
which is not given locally by its Taylor series
around any point; we present a $C^1_\C$-map
$f\!: \C\to E$ to a metrizable, complete,
non-locally convex topological vector space
which is not $C^2_\C$; and we present
a non-zero, compactly supported $C^\infty_\C$-map
$\C\to E$, the existence of which
demonstrates that the Identity Theorem
fails for suitable $C^\infty_\C$-maps
into non-locally convex spaces.\\[3mm]
Mappings
between open subsets of the field $\Q_p$
of $p$-adic numbers with similar pathological
properties
are known in non-archimedian analysis
(see \cite{Sch}).
The author also drew inspiration from
\cite[Ex.\,II.2.7]{NNa},
where
it is shown
that the map $[0,1]\to L^{1/2}[0,1]$, $t\mto \one_{[0,t]}$
is differentiable at each point (in the ordinary sense),
with vanishing derivative.\\
$\;$\\
{\bf Differential calculus.}
We shall not give the general definition
of $C^r_\K$-maps between open subsets of
topological $\K$-vector spaces here.
Rather, we
shall work with a simpler definition for the special case
of curves (cf.\ also \cite{Sch}),
which is equivalent
to the general definition (cf.\ \cite[Prop.\,6.9]{Ber}).
All topological vector spaces are assumed Hausdorff.\\[3mm]
{\bf Definition.\/}
Let $\K$ be a non-discrete topological field,
and $f\!: U\to E$ be a mapping
from an open subset $U$ of $\K$
to a topological $\K$-vector space~$E$.
The map $f$ is said to be {\em of class $C^0_\K$\/}
if it is continuous; in this case, define $f^{<0>}:=f$.
We call $f$ {\em of class $C^1_\K$\/}
if it is continuous and if there exists
a continuous map $f^{<1>}\!: U\times U\to E$
such that
\[
f^{<1>}(x_1,x_2)=\frac{1}{x_1-x_2}(f(x_1)-f(x_2))
\]
for all $x_1,x_2\in U$ such that $x_1\not=x_2$.
Recursively, having defined mappings of class
$C^j_\K$ and associated maps $f^{<j>}\!: U^{j+1}\to E$
for $j=0,\ldots, k-1$ for some $k\in \N$,
we call $f$ {\em of class $C^k_\K$\/}
if it is of class $C^{k-1}_\K$
and if there exists a continuous map
$f^{<k>}\!: U^{k+1}\to E$
such that
\[
f^{<k>}(x_1,x_2,\ldots, x_{k+1})=
\frac{1}{x_1-x_2}\left(f^{<k-1>}(x_1,x_3,\ldots, x_{k+1})
-f^{<k-1>}(x_2,x_3,\ldots, x_{k+1})\right)
\]
for all $x_1,\ldots,x_{k+1}\in U^{k+1}$
such that $x_1\not=x_2$.
The map $f$ is {\em of class $C^\infty_\K$\/}
(or {\em smooth\/})
if it is of class $C^k_\K$
for all $k\in \N_0$.\\[3mm]
Here $f^{<k>}$ is uniquely determined,
and $f^{<k>}$ is symmetric in its $k+1$
variables.
Furthermore, $k!\, f^{<k>}(x,\ldots,x)=\frac{d^kf}{dx^k}(x)=:f^{(k)}(x)$,
for all $x\in U$ (cf.\ \cite{Ber}).\\[3mm]
\noindent
{\bf Example~1.\/}
Let
$d\nu(x)=\frac{1}{\sqrt{\pi}}\, e^{-x^2}\, dx$ be the Gauss
measure on $\R$,
and
$\mu:=\nu\tensor \nu$
be the Gauss measure
on~$\C=\R^2$.
Let $E:=L^0(\C,\mu)$ be the complex topological
vector space
of equivalence classes of measurable
complex-valued functions on~$\C$ (modulo functions vanishing
almost everywhere), equipped with
the topology of convergence in measure.
Thus, a basis of open zero-neighbourhoods of~$E$
is given by the sets $W_k$ for $k\in \N$,
where $W_k$
consists of those equivalence classes
of measurable maps $\gamma\!: \C\to\C$
such that
$\mu\left(\{z\in \C\!: |\gamma(z)|\geq \frac{1}{k}\}\right)<\frac{1}{k}$.
It is well-known that $E$ is a metrizable,
complete, non-locally convex  topological vector space,
which does not admit any non-zero continuous
linear functionals (cf.\ \cite{KPR}).
Given $\gamma \in L^0(\C,\mu)$ and a closed subset $A\sub \C$,
we say that $\gamma$ is {\em supported in~$A$\/}
if $\gamma$ vanishes $\mu$-almost everywhere
on the complement of~$A$. The {\em support of~$\gamma$\/}
is the smallest closed set in which $\gamma$ is supported.
In the following, $\one_A\!: \C\to \{0,1\}$
denotes the characteristic function of a measurable
subset $A\sub \C$.
\\[3mm]
Consider the mapping
\[
f\!: \C\to E, \quad f(z):=\one_{A(z)},
\]
where $A(z):=\{w\in \C\!:
\,\mbox{$\Repart(w)\leq \Repart(z)$ and $\Impart(w)\leq \Impart(z)$}\,\}$.
Then $f$ has the following properties:
\begin{prop}\label{maj2}
$f\!: \C\to E$ is of class $C^\infty_\C$,
injective,
and $f^{(j)}(z)=0$ for all
$j\in \N$ and $z\in \C$.
In particular, $f$ is not given locally by
its Taylor series around any point.
\end{prop}
\begin{proof}
Apparently $f$ is injective. If
we can show that $f$ is $C^\infty_\C$,
with $f^{(j)}=0$ for all $j\in \N_0$,
then, given any $z_0\in \C$,
the Taylor series of~$f$ at~$z_0$
will only consist of the 0th order term,
and hence describes the
function which is constantly $f(z_0)$.
It therefore does not coincide with the injective function~$f$
on any neighbourhood of~$z_0$.\\[3mm]
Thus, to complete the proof of the proposition,
it suffices to establish the following assertions,
by induction on $k\in\N_0$:
\begin{itemize}
\item[\n (a)]
$f$ is of class $C^k_\C$;
\item[\n (b)]
For all $j\in \N$ such that $j\leq k$, we have $f^{(j)}=0$,
and, for all $z_1,\ldots, z_{j+1}\in \C$,
the element
$f^{<j>}(z_1,\ldots, z_{j+1})\in E=L^0(\C,\mu)$ is supported in
\[
([x_*,x^*]+i\R) \;\cup \; (\R+i[y_*,y^*])\,,\quad \mbox{where}
\]
\begin{equation}\label{xstar}
\begin{array}{cccccc}
x_* & := & \min\{\Repart(z_1),\ldots,\Repart(z_{j+1})\}, &
x^*& := & \max\{\Repart(z_1),\ldots,\Repart(z_{j+1})\},\\
y_*&:= &\min\{\Impart(z_1),\ldots,\Impart(z_{j+1})\},&
y^* & := & \max\{\Impart(z_1),\ldots,\Impart(z_{j+1})\}\,.
\end{array}
\end{equation}
\end{itemize}
{\em The case $k=0$.} Given $z_1\in \C$,
let us show that $f$ is continuous at~$z_1$.
To this end, let $z_2\in \C$; define $x_*,x^*,y_*,y^*$
as in (\ref{xstar}) (taking $j:=1$).
Then the symmetric difference $A(z_1)\oplus A(z_2):=
(A(z_1)\setminus A(z_2))\,\cup
\, (A(z_2)\setminus A(z_1))$
of the sets $A(z_1)$ and $A(z_2)$ is a subset
of $([x_*,x^*]+i\R) \,\cup \, (\R+i[y_*,y^*])$,
whence
\[
\mu(A(z_1)\oplus A(z_2))\leq |x^*-x_*|+|y^*-y_*|\leq 2\,|z_2-z_1|
\]
(using Fubini's Theorem). Note that $f(z_2)-f(z_1)$
is supported in $A(z_1)\oplus A(z_2)$.
As the measure of this set tends to~$0$
as $z_2\to z_1$, we see that $f(z_2)\to f(z_1)$
in~$E$. Thus $f$ is continuous.\\[3mm]
{\em Induction step.}
Let $k\in \N$, and suppose that (a) and (b) hold
when $k$ is replaced with $k-1$.
Then $f$ is of class $C^{k-1}_\C$.
In order that $f$ be $C^k_\C$,
with $f^{(k)}=0$, in view of
La.\,10.5, La.\,10.7 and
Prop.\,6.2 in~\cite{Ber},
we only need to show that
\[
g_n:=\frac{1}{z_{n,1}-z_{n,2}}
\left( f^{<k-1>}(z_{n,1},z_{n,3},\ldots,z_{n,k+1})
-f^{<k-1>}(z_{n,2},z_{n,3},\ldots,z_{n,k+1})\right)\to 0
\]
in $E$
as $n\to\infty$, for every
sequence $(z_n)_{n\in \N}$
of elements $z_n=(z_{n,1},\ldots,z_{n,k+1})\in \C^{k+1}$
which converges to a diagonal element $(z,z,\ldots, z)\in \C^{k+1}$
for some $z\in \C$,
where $z_{n,a}\not=z_{n,b}$
whenever $a\not=b$.
Given $n\in \N$,
define $x_{n,*}$, $x_n^*$, $y_{n,*}$ and
$y_n^*$ along the lines of (\ref{xstar}),
using the elements $z_{n,1},\ldots, z_{n,k+1}$
(thus $j=k$).
Note that, as a consequence of (b) for $k$ replaced with
$k-1$
(valid by the induction
hypothesis),
each of the elements
$f^{<k-1>}(z_{n,2},z_{n,3},\ldots,z_{n,k+1})$,
$f^{<k-1>}(z_{n,1},z_{n,3},\ldots,z_{n,k+1})\in E$
is supported in
\[
B_n:=([x_{n,*}, x_n^*]+i\R)\; \cup\;(\R+i[y_{n,*},y_n^*])\, .
\]
Hence also $g_n$ is supported in $B_n$.
Since
\[
\mu(B_n)\leq 4\max\{|z_{n,1}-z|,\ldots,|z_{n,k+1}-z|\}\to 0\quad
\mbox{as $\;n\to\infty$,}
\]
we deduce that $\lim_{n\to\infty}g_n=0$ in~$E$,
as required. Thus $f$ is $C^k_\C$,
and $k!\, f^{<k>}(z,\ldots,z)= f^{(k)}(z)=0$
for all $z\in \C$.\\[3mm]
It only remains to prove the assertion concerning the
supports. To this end, let $z_1,\ldots, z_{k+1}$
$\in \C$.
If all of $z_1,\ldots, z_{k+1}$ coincide,
then $f^{<k>}(z_1,\ldots,z_{k+1})=0$
by the preceding, and this is an element with empty
support, which therefore is contained
in the desired set.
Now suppose that
$z_a\not= z_b$ for some $a,b$.
By symmetry of $f^{<k>}$ in its $k+1$
variables, we may assume that $z_1\not=z_2$.
Then
\[
f^{<k>}(z_1,\ldots, z_{k+1})=\frac{1}{z_1-z_2}
\left(f^{<k-1>}(z_1,z_3,\ldots, z_{k+1})-
f^{<k-1>}(z_2,z_3,\ldots, z_{k+1})\right)
\]
and, as in the preceding part of the proof,
we deduce from the induction hypothesis
that this element is supported in the
desired set.
\end{proof}
\begin{cor}\label{cor1}
Consider $E$ as a real topological vector space.
Then $g:=f|_\R\!: \R\to E$
is an injective
$C^\infty_\R$-curve
whose derivative $g'$ vanishes identically.\Punkt
\end{cor}
{\bf Example~2.\/} We retain $\mu$ and $E$ as in Example~1,
but consider now the mapping
\[
f\!: \C\to E, \quad f(z):=\one_{A(z)},
\]
where $A(z):=\{w\in \C\!: |z|\leq |w|\leq 1\}$.
Then $f$ has the following properties:
\begin{prop}\label{maj3}
$f\!: \C\to E$ is
of class $C^\infty_\C$, non-zero,
and $f$ has compact support.
\end{prop}
\begin{proof}
Clearly $f(z)=0$ for all $z\in \C$ such that $|z|\geq 1$,
entailing that $f$ is compactly supported.
Furthermore, $f\not=0$.
Given real numbers $0\leq r\leq R$,
let
\[
K(r,R):=\{w\in \C\!: r\leq |w|\leq R\}
\]
be the closed annulus with inner radius~$r$
and outer radius~$R$ in~$\C$.
Then
\begin{equation}\label{volannul}
\mu(K(r,R))\leq R-r\,.
\end{equation}
Indeed, we have
\begin{eqnarray*}
\mu(K(r,R))&= &\int_r^R\int_0^{2\pi}{\textstyle \frac{s}{\pi}}\,
e^{-s^2}\,
d\phi \,ds=e^{-r^2}-e^{-R^2}\\
&=&(r-R)\cdot(-2\xi \,e^{-\xi^2}) =(R-r)\,2\xi\,e^{-\xi^2}\leq R-r
\end{eqnarray*}
for some $\xi\in [r,R]$,
using the Mean Value Theorem to pass
to the second line, and using that $2te^{-t^2}\leq \sqrt{\frac{2}{e}}<1$
for all $t\in [0,\infty[$, by an elementary calculation.\\[3mm]
The assertion of the proposition will follow if we can prove
the following claims, by induction on $k\in \N_0$:
\begin{itemize}
\item[\n (a)]
$f$ is of class $C^k_\C$;
\item[\n (b)]
For all $j\in \N$ such that $j\leq k$, the map $f^{(j)}$
vanishes,
and
$f^{<j>}(z_1,\ldots, z_{j+1})\in E$ is supported in
the annulus
$K(r_*,r^*)$, where
\[
r_* := \min\{|z_1|,\ldots,|z_{j+1}|\}, \quad\quad
r^*:= \max\{|z_1|,\ldots,|z_{j+1}|\}\,,
\]
for all $z_1,\ldots, z_{j+1}\in \C$.
\end{itemize}
In view of (\ref{volannul}),
an apparent adaptation of
the proof of Prop.\,\ref{maj2}
establishes these claims.
\end{proof}
{\bf Remark.\/}
Proposition~\ref{maj3} entails
that the identity theorem for analytic mappings
becomes invalid when analytic maps are replaced
with $C^\infty_\C$-maps into
non-locally convex spaces.\\[3mm]
{\bf Example~3.} We retain $\nu$ and $\mu$ as in Example~1
but consider
the complex topological vector space
$E:=L^p(\C,\mu)$ of equivalence classes
of complex-valued $L^p$-functions on $\C$ now,
where $p\in \;]\frac{1}{2},1[$.
Consider
\[
f\!: \C\to E, \quad f(z):=\one_{A(z)},
\]
where $A(z):=\{w\in \C\!: \Repart(w)\leq \Repart(z)\}$.
The map $f$ has the following properties:
\begin{prop}\label{maj4}
$f\!: \C\to E$ is
of class $C^1_\C$ and $f'$ vanishes,
but $f$ is not of class $C^2_\C$.
\end{prop}
\begin{proof}
Given real numbers $a\leq b$, define
$S(a,b):=\{w\in \C\!: a<\Repart(w)\leq b\}$.
Then
\begin{equation}\label{muS}
\mu(S(a,b))=\nu(]a,b])\leq b-a\,.
\end{equation}
Given $z_1,z_2\in \C$, where $\Repart(z_1)\leq \Repart(z_2)$
without loss of generality,
we have
\[
f(z_2)-f(z_1)=\one_{S(\repart(z_1),\repart(z_2))}\,,
\]
where $\mu(S(\Repart(z_1),
\Repart(z_2)))\leq \Repart(z_2)-\Repart(z_1)
\leq |z_2-z_1|$.
We easily deduce that $f$ is continuous.
Assuming that $z_1\not=z_2$ here, we have
\[
\int_\C
\left|\frac{f(z_2)(w)-f(z_1)(w)}{z_2-z_1}\right|^p\, d\mu(w)
=|z_2-z_1|^{-p}\cdot \mu(S(\Repart(z_1),\Repart(z_2)))
\leq |z_2-z_1|^{1-p}\to 0
\]
as $|z_2-z_1|\to 0$, showing that $\frac{1}{z_2-z_1}(f(z_2)-f(z_1))\to 0$
in~$E$ whenever $|z_2-z_1|\to 0$.
Thus $f^{<1>}(z_1,z_2):=0$ if $z_1=z_2$,
$f^{<1>}(z_1,z_2):=\frac{1}{z_2-z_1}(f(z_2)-f(z_1))$
if $z_1\not=z_2$
defines a continuous function $f^{<1>}\!: \C\to E$,
showing that $f$ is $C^1_\C$ with $f'(z)=f^{<1>}(z,z)=0$
for all $z\in \C$.\\[3mm]
Let $c:=\frac{1}{e\sqrt{\pi}}$;
then $\frac{1}{\sqrt{\pi}}e^{-t^2}\geq c$ for all
$t\in [0,1]$.
For $t\in [0,\frac{1}{2}]$, we have
\[
{\textstyle\frac{1}{t}}\left( f^{<1>}(t,2t)-f^{<1>}(0,2t)\right)
={\textstyle\frac{1}{t^2}}\left(\one_{S(t,2t)}-{\textstyle\frac{1}{2}}
\one_{S(0,2t)}\right)\,,
\quad\mbox{whence}
\]
\begin{eqnarray*}
\lefteqn{\!\!\!\!\!\!\!
\int_\C \left|{\textstyle\frac{1}{t}}\left( f^{<1>}(t,2t)-f^{<1>}(0,2t)\right)
(w)\right|^p\,
d\mu(w)\quad}\\
\quad &=&{\textstyle \left(\frac{1}{2t^2}\right)^p}\cdot \mu(S(0,2t))
={\textstyle \left(\frac{1}{2t^2}\right)^p}\cdot \nu(]0,2t])
\geq
{\textstyle \left(\frac{1}{2t^2}\right)^p} 2t\, c
=2^{1-p}t^{1-2p}\, c\to \infty
\end{eqnarray*}
as $t\to 0$.
The map $E\to [0,\infty[$, $\gamma\mto \int_\C |\gamma(w)|^p\, d\mu(w)$
being continuous, we deduce that
$\lim_{t\to 0}
\frac{1}{t}\left( f^{<1>}(t,2t)-f^{<1>}(0,2t)\right)$
cannot exist in~$E$. Thus $f$ is not $C^2_\C$.
\end{proof}
\begin{cor}\label{cor2}
Consider $E$ as a real topological vector space.
Then $g:=f|_\R\!: \R\to E$
is a
$C^1_\R$-curve
whose derivative $g'$ vanishes identically.
However, $g$ is not $C^2_\R$.\Punkt
\end{cor}
\noindent
{\footnotesize
{\bf Helge Gl\"{o}ckner}, TU~Darmstadt, FB~Mathematik~AG~5,
Schlossgartenstr.\,7, 64289 Darmstadt, Germany.\\
E-Mail: gloeckner@mathematik.tu-darmstadt.de}
\end{document}